\documentclass[11pt]{article}

\usepackage{amsmath,amsbsy,amssymb,latexsym,amsthm,dsfont}
\usepackage{a4wide}

\usepackage{cite}
\usepackage{url}

\newtheorem{Theorem}{Theorem}
\newtheorem{Definition}[Theorem]{Definition}
\newtheorem{Remark}[Theorem]{Remark}

\newenvironment{keywords}{\begin{center}
\begin{minipage}[c]{13.4cm} {\bf Keywords:}} {\end{minipage}
\end{center}}

\newenvironment{msc}{\begin{center}
\begin{minipage}[c]{13.4cm} {\bf MSC 2010:}} {\end{minipage}
\end{center}}

\newenvironment{pacs}{\begin{center}
\begin{minipage}[c]{13.4cm} {\bf PACS:}} {\end{minipage}
\end{center}}

\begin{document}

\title{Nondifferentiable variational principles\\
in terms of a quantum operator\thanks{Submitted 24-Apr-2011;
revised 18-Jun-2011; accepted 20-Jun-2011;
for publication in \emph{Mathematical Methods in the Applied Sciences}.}}

\author{Ricardo Almeida\\
\texttt{ricardo.almeida@ua.pt}
\and Delfim F. M. Torres\\
\texttt{delfim@ua.pt}}

\date{Department of Mathematics,
University of Aveiro,
3810-193 Aveiro, Portugal
}

\maketitle


\begin{abstract}
We develop Cresson's nondifferentiable calculus of variations
on the space of H\"{o}lder functions.
Several quantum variational problems are considered: with and without constraints,
with one and more than one independent variable, of first and higher-order type.
\end{abstract}


\begin{msc}
26A27; 26B20; 39A13; 49K05; 49K10; 49S05.
\end{msc}



\section{Introduction}

Many physical phenomena are described by nondifferentiable functions.
For instance, generic trajectories of quantum mechanics
are nondifferentiable \cite{Feynman}.
An important issue consists to find stationary conditions for integral functionals defined
on sets of functions that are not necessarily differentiable in the classical
sense. Several different approaches to deal with nondifferentiability
are being followed in the literature, including
the time scale approach, which typically deal with delta or nabla differentiable
functions \cite{AlmeidaNabla,Bartos,MyID:159,malina,MyID:142,NataliaHigherOrderNabla,MyID:198},
the fractional approach, allowing to consider
functions that have no first order derivative
but have fractional derivatives of all orders less than one
\cite{MyID:182,RicAML,MR2557007,El-Nabulsi,El-Nabulsi:Torres:JMP,withRui:hZ:JDDE,gasta:07},
and the quantum approach, particularly useful to model physical systems
\cite{Bang04,Bang05,withMiguel01,Cresson,Torres,CressonGreff,Cresson1,Jannussis,Malinowska}.

Quantum derivatives play a leading role in the understanding of complex
physical systems. In 1992 Nottale introduced the theory of scale-relativity
without the hypothesis of space-time differentiability
\cite{Nottale:1992,Nottale1,Nottale2}. In \cite{Cresson} Cresson presents a quantum calculus
defined on a set of H\"{o}lder functions based on the
$h$-operator
$$\Delta_{h} f(x)=\frac{f(x+h)-f(x)}{h}.$$
Main result of \cite{Cresson} gives a nondifferentiable Euler--Lagrange equation
for the basic problem of the calculus of variations.
Cresson's calculus of variations has been further developed in
\cite{Torres}, where Lagrangian and Hamiltonian versions
of Noether's theorem are proved, and in \cite{Almeida1}
where the authors study isoperimetric problems. More recently,
Cresson and Greff \cite{CressonGreff,Cresson1}
improved the previous approach of \cite{Cresson}.

Our paper is organized as follows. In Section~\ref{sec:def} we recall
the definition of quantum derivative (Definition~\ref{def:ourHD})
and its main properties are reviewed, including Leibniz's and Barrow's formulas
(Theorems~\ref{LeibnizRule} and \ref{Barrow}, respectively).
We also introduce the concept of higher-order and partial
quantum derivatives. In Section~\ref{sec:optimization1} we consider
integral functionals with Lagrangians containing
the  quantum derivative, defined on a set of H\"{o}lder functions.
We present a necessary and sufficient extremality condition of Euler--Lagrange type
(Theorem~\ref{ELtheorem1}) as well as natural boundary conditions
(Theorem~\ref{thm:sec3:nbc}). In Section~\ref{sec:iso} we exhibit a generalization
of the  quantum Euler--Lagrange equation when we are in presence of an integral constraint
of the same type as the quantum action functional (Theorem~\ref{thm:EL:ISO:P}).
In Section~\ref{sec:parameter} we consider dependence of the action
on a complex parameter (Theorem~\ref{thm:EL:compParameter}).
The case where the Lagrangian contains higher-order quantum derivatives
is considered in Section~\ref{sec:higherEL}, and a higher-order
Euler--Lagrange type equation is proved (Theorem~\ref{higherEL}).
The rest of the paper considers a generalization for admissible
functions of two independent variables. We begin by proving
in Section~\ref{sec:green} a nondifferentiable quantum version of Green's theorem
(Theorem~\ref{thm:HQGT}) and a quantum analogous of integration by parts for
double integrals (formula \eqref{eq:intParts2d}),
from which we deduce in Section~\ref{sec:optimization2}
a necessary and sufficient condition for a function
to be a quantum extremal for double integrals
(Theorem~\ref{ELtheorem}) and respective natural
boundary conditions (Theorem~\ref{ELtheorem:di}).
We end the paper with Section~\ref{sec:example}, considering
a membrane with potential energy given by a variational double integral
depending on  quantum partial derivatives.
The equilibrium of the membrane is given
as the  quantum Euler--Lagrange extremal.


\section{Definitions and basic formulas}
\label{sec:def}

We begin by reviewing the quantum calculus as in \cite{CressonGreff,Cresson1}.
Along the text $\alpha, \beta \in(0,1)$,
$h \in(0,1)$ with $h \ll 1$ and $\sigma=\pm 1$.
For a continuous function $f:I\subseteq \mathbb{R} \to \mathbb{R}$,
the $h$-derivative of $f$ is defined by the quotient
$${D^\sigma_{h}f(x)}=\sigma\frac{f(x+\sigma h)-f(x)}{h}$$
whenever it is defined.
Observe that, in case of considering differentiable functions,
we can obtain the standard derivative taking the limit
$$\lim_{h\to 0}D^\sigma_{h}f(x)=f'(x).$$

\begin{Definition}
The quantum derivative of $f$, with respect to $h$, is defined by
\begin{equation}
\label{eq:def:cp}
\frac{{\Box_{h}}f}{\Box x}(x)
=\frac12 \left[ \left( {D^{+1}_{h}f(x)}+{D^{-1}_{h}f(x)}\right)
+i\left( {D^{+1}_{h}f(x)} - {D^{-1}_{h}f(x)} \right) \right].
\end{equation}
\end{Definition}

As before, for differentiable functions we obtain $f'(x)$ from
\eqref{eq:def:cp} as $h\to0$.
For complex valued functions we put
$$
\frac{{\Box_{h}}f}{\Box x}(x)
= \frac{{\Box_{h}}\mbox{Re}f}{\Box x}(x)
+ i \frac{{\Box_{h}}\mbox{Im}f}{\Box x}(x).
$$

Let ${C^0_{conv}}(I\times (0,1),\mathbb{R})\subseteq {C^0}(I\times (0,1),\mathbb{R})$
be the set of functions for which the limit
$$\lim_{h\to 0}f(x,h)$$
exists for any $x\in I$, and let $E$
be a complementary of ${C^0_{conv}}(I\times (0,1),\mathbb{R})$
in ${C^0}(I\times (0,1),\mathbb{R})$.
Define the projection map $\pi$ by
$$
\begin{array}{lcll}
\pi: & {C^0_{conv}}(I\times (0,1),\mathbb{R})
\oplus E & \to & {C^0_{conv}}(I\times (0,1),\mathbb{R})\\
& f_{conv}+f_E  & \mapsto & f_{conv}
\end{array}
$$
and the operator $\left< \cdot \right>$ by
$$
\begin{array}{lcll}
\left< \cdot \right>: & {C^0}(I\times (0,1),\mathbb{R}) & \to & {C^0}(I,\mathbb{R})\\
& f & \mapsto & \left< f \right>: x\mapsto
\displaystyle\lim_{h\to 0}\pi(f)(x,h).
\end{array}
$$

Finally, we arrive to the main concept introduced in \cite{CressonGreff}:
the quantum derivative of $f$ (without the dependence of $h$).

\begin{Definition}
\label{def:ourHD}
The quantum derivative of $f$ is defined by the rule
\begin{equation}
\label{eq:scaleDer}
\frac{\Box f}{\Box x}=\left< \frac{{\Box_{h}}f}{\Box x} \right>.
\end{equation}
\end{Definition}

The scale derivative \eqref{eq:scaleDer} has some nice properties.
Namely, it satisfies a Leibniz and a Barrow rule.
First let us recall the definition of H\"{o}lderian function.

\begin{Definition}
Let $f\in C^0(I,\mathbb{R})$. We say that $f$ is H\"{o}lderian of H\"{o}lder exponent
$\alpha$ if there exists a constant $C>0$ such that,
for all $x$ and $s$, the inequality
$$
|f(x+s)-f(x)|\leq C |s|^\alpha
$$
holds. The set of  H\"{o}lderian functions of H\"{o}lder
exponent $\alpha$ is denoted by $H^\alpha(I, \mathbb{R})$.
\end{Definition}

\begin{Theorem}[The quantum Leibniz rule \cite{CressonGreff}]
\label{LeibnizRule}
For $f\in H^\alpha(I, \mathbb{R})$ and $g\in H^\beta(I, \mathbb{R})$,
$\alpha+\beta>1$, one has
\begin{equation*}
\frac{\Box (f.g)}{\Box x}(x)
=\frac{\Box f}{\Box x}(x).g(x)+f(x).\frac{\Box g}{\Box x}(x).
\end{equation*}
\end{Theorem}

\begin{Theorem}[The quantum Barrow rule \cite{CressonGreff}]
\label{Barrow}
Let $f\in C^0([a,b],\mathbb{R})$ be such that $\Box f / \Box x$ is continuous, and
\begin{equation}
\label{nec_condition}
\lim_{h\to0} \int_a^b \left(\frac{\Box_h f}{\Box x}\right)_E(x)dx=0.
\end{equation}
Then,
$$
\int_a^b \frac{\Box f}{\Box x}(x)\, dx=f(b)-f(a).
$$
\end{Theorem}

We can generalize the previous notion to include higher-order
and partial derivatives. An higher-order
quantum derivative is defined as in the standard case:
given $n\in \mathbb{N}$, the  quantum derivative of order
$n$ of $f$, $\frac{\Box^n f}{\Box x^n}$, is defined recursively by
\begin{equation*}
\begin{split}
\frac{\Box^1 f}{\Box x^1}&= \frac{\Box f}{\Box x},\\
\frac{\Box^n f}{\Box x^n}&=\frac{\Box }{\Box x}\left(\frac{\Box^{n-1} f}{\Box x^{n-1}}\right),
\quad n = 2, 3, \ldots
\end{split}
\end{equation*}
For convenience of notation, we use the convention
$\frac{\Box^0 f}{\Box x^0}= f$. We note that
if $f$ is a function of class $C^n$, then
$\frac{\Box^n f}{\Box x^n} = \frac{d^n f}{d x^n} = f^{(n)}$.

We now consider the two dimensional case. The objective is to introduce the concept
of partial  quantum derivative. All the previous definitions and results
can be easily adapted to this case, doing the necessary adjustments.
We proceed doing a brief exposition on the subject.
Let $\Omega\subseteq\mathbb{R}^2$. For $\sigma=\pm 1$ we define the quantities
$$
\frac{{\partial^\sigma_{h}}f}{\partial x_1}(x_1,x_2)
=\sigma\frac{f(x_1+\sigma h,x_2)-f(x_1,x_2)}{h}
$$
and
$$
\frac{{\partial^\sigma_{h}}f}{\partial x_2}(x_1,x_2)
=\sigma\frac{f(x_1,x_2+\sigma h)-f(x_1,x_2)}{h}.
$$
If $f$ is a differentiable function, then
we obtain the usual partial derivatives
$\partial f/\partial x_1$ and $\partial f/\partial x_2$ as $h\to 0$.
The partial  $h$-quantum derivatives are defined by
$$
\frac{{\Box_{h}}f}{\Box x_j}(x_1,x_2)
=\frac12 \left[ \left( \frac{{\partial^{+1}_{h}}f}{\partial x_j}(x_1,x_2)
+\frac{{\partial^{-1}_{h}}f}{\partial x_j}(x_1,x_2)\right)
+i\left( \frac{{\partial^{+1}_{h}}f}{\partial x_j}(x_1,x_2)
-\frac{{\partial^{-1}_{h}}f}{\partial x_j}(x_1,x_2)\right) \right],
$$
$j\in\{1,2\}$.

Given a function $f(\cdot,\cdot)$, we say that $f$ belongs to
$H^\alpha(\Omega,\mathbb{R})$ if both $f(\cdot,x_2)$ and $f(x_1,\cdot)$
are H\"{o}lderian of H\"{o}lder exponent $\alpha$ for all $x_1$ and $x_2$.

The definition of the partial  quantum derivatives
$\Box f/\Box x_j$ of $f$, $j=1,2$, is clear and is left to the reader.


\section{Extremals for functions with one independent variable:
the  quantum Euler--Lagrange equation}
\label{sec:optimization1}

In the classical context of the calculus of variations,
for example in classical mechanics, one studies necessary
and sufficient optimality conditions for functionals/actions defined
on a set of differentiable functions $y$, where the Lagrangian $L$
depends on the independent variable $x$, function $y$, and the derivative $y'$ of $y$.
In our context we deal with nondifferentiable functions $y$,
and we replace $y' = d y/d x$ by the  quantum derivative
$\Box y/\Box x$. We consider the following functional:
\begin{equation}
\label{funct}
\begin{array}{lcll}
\Phi: & \partial H^\alpha(I,\mathbb{R}) & \to & \mathbb{C}\\
& y & \mapsto & \displaystyle \int_a^b L\left(x,y(x),\frac{\Box y}{\Box x}(x)\right)\, dx,
\end{array}
\end{equation}
where $\partial H^\alpha(I,\mathbb{R})$ denotes the subspace of
$H^\alpha(I,\mathbb{R})$ for which $y(a)$ and $y(b)$ take given
fixed values $y_a$ and $y_b$, \textrm{i.e.}, $y(a)=y_a$ and $y(b)=y_b$.
We assume that the Lagrangian
$L=L(x,y,v):[a,b] \times \mathbb{R}\times \mathbb{C} \to \mathbb{C}$
is differentiable in all its arguments with continuous partial derivatives:
$L(\cdot,\cdot,\cdot) \in C^1$.

\begin{Remark}
If in \eqref{funct} we replace the domain of $\Phi$ by $C^1(I,\mathbb{R})$,
then we obtain the usual functional of variational calculus:
$$
\Phi(y)=\displaystyle \int_a^b L(x,y(x),y'(x))\, dx.
$$
\end{Remark}

Let $\alpha, \beta \in(0,1)$ be such that $\alpha+\beta>1$ and $\beta \geq \alpha$.
A variation of $y\in \partial H^\alpha(I,\mathbb{R})$
is another function of $\partial H^\alpha(I,\mathbb{R})$
of the form $y+ \varepsilon w$, with $w \in H^\beta(I, \mathbb{R})$ and
$\varepsilon$ a small real number.

\begin{Definition}
We say that $y$ is a quantum extremal of \eqref{funct}
if for any $w\in H^\beta(I, \mathbb{R})$
$$
\frac{d}{d\varepsilon}\left.\Phi(y + \varepsilon w)\right|_{\varepsilon = 0}=0.
$$
\end{Definition}

For convenience of notation, we introduce the operator $[\cdot]$ defined by
$$[y](x) =\left(x,y(x),\frac{\Box y}{\Box x}(x)\right).$$

\begin{Theorem}[The quantum Euler--Lagrange equation \cite{CressonGreff}]
\label{ELtheorem1}
Let $\alpha, \beta \in(0,1)$ be such that $\alpha+\beta>1$ and $\beta \geq \alpha$.
If $y \in \partial H^\alpha(I,\mathbb{R})$ satisfies
$\Box y/\Box x \in H^\alpha(I,\mathbb{R})$ and
$$\frac{\partial L}{\partial v}[y](x)\cdot w(x)$$
satisfies condition \eqref{nec_condition} for all $w \in H^\beta(I, \mathbb{R})$,
then function $y$ is a quantum extremal of $\Phi$ on $H^\beta(I, \mathbb{R})$ if and only if
\begin{equation}
\label{ELequation1}
\frac{\partial L}{\partial y}[y](x)
-\frac{\Box}{\Box x}\left(\frac{\partial L}{\partial v}[y]\right)(x)=0
\end{equation}
for all $x\in [a,b]$.
\end{Theorem}
From now on we will be always assuming that
$$
\frac{\partial L}{\partial v}[y](x)\cdot w(x)
$$
satisfies \eqref{nec_condition} for all $w \in H^\beta(I, \mathbb{R})$.

Another important problem is to find the extremals of a certain functional
$\Psi$ in the case we have no boundary constraints on the set of admissible functions
that define the domain of $\Psi$. As we show below, in this case
besides the  quantum Euler--Lagrange equation two more conditions,
called the quantum natural boundary conditions, need to be satisfied.
Let $\Psi$ be the functional defined by
$$
\begin{array}{lcll}
\Psi: & H^\alpha(I,\mathbb{R}) & \to & \mathbb{C}\\
& y & \mapsto & \Phi(y)
= \displaystyle \int_a^b L\left(x,y(x),\frac{\Box y}{\Box x}(x)\right)\, dx.
\end{array}
$$

\begin{Theorem}[The quantum Euler--Lagrange and natural boundary conditions]
\label{thm:sec3:nbc}
Let $y \in H^\alpha(I,\mathbb{R})$ be such that
$\Box y/\Box x \in H^\alpha(I,\mathbb{R})$.
Function $y$ is a quantum extremal of $\Psi$ on $H^\beta(I, \mathbb{R})$
if and only if the following conditions hold:
\begin{enumerate}
\item $\displaystyle\frac{\partial L}{\partial y}[y](x)
-\frac{\Box}{\Box x}\left(\frac{\partial L}{\partial v}[y]\right)(x)=0$ for all $x\in [a,b]$;
\item $\displaystyle\frac{\partial L}{\partial v}\left( a,y(a),\frac{\Box y}{\Box x}(a) \right)=0$;
\item $\displaystyle\frac{\partial L}{\partial v}\left( b,y(b),\frac{\Box y}{\Box x}(b)\right)=0$.
\end{enumerate}
\end{Theorem}

\begin{proof}
Similarly to the proof of Theorem~\ref{ELtheorem1}, we begin by considering
a variation of $y$ of the type $y+ \varepsilon w$, with $w\in H^\beta(I, \mathbb{R})$.
Using Theorems~\ref{LeibnizRule} and \ref{Barrow}, it follows that for all $w\in H^\beta(I, \mathbb{R})$,
\begin{equation*}
\begin{split}
0 & = \frac{d}{d\varepsilon}\left.\Phi(y + \varepsilon w)\right|_{\varepsilon = 0}\\
& =\displaystyle \int_a^b \left[ \frac{\partial L}{\partial y}[y](x)\cdot w(x)
+\frac{\partial L}{\partial v}[y](x)\cdot \frac{\Box w}{\Box x}(x) \right] dx\\
 &=\displaystyle \int_a^b \left[ \frac{\partial L}{\partial y}
-\frac{\Box}{\Box x}\left(\frac{\partial L}{\partial v} \right) \right]
\cdot w \, dx+\int_a^b\frac{\Box}{\Box x}\left(\frac{\partial L}{\partial v}\cdot w\right) \, dx\\
&=\displaystyle \int_a^b \left[ \frac{\partial L}{\partial y}
-\frac{\Box}{\Box x}\left(\frac{\partial L}{\partial v} \right) \right] \cdot w \, dx
+ \frac{\partial L}{\partial v}\left(b,y(b),\frac{\Box y}{\Box x}(b)\right)\cdot w(b)
- \frac{\partial L}{\partial v}\left(a,y(a),\frac{\Box y}{\Box x}(a)\right)\cdot w(a).
\end{split}
\end{equation*}
Since the last equality is valid for any $w$, we can choose those such that
$w(a)=w(b)=0$, obtaining the Euler--Lagrange equation \eqref{ELequation1}
given in item 1. Moreover, choosing functions $w$ for which
$w(a)=0$ and $w(b)\not=0$ or $w(b)=0$ and $w(a)\not=0$,
we deduce the remaining two conditions.
\end{proof}


\section{The  quantum isoperimetric problem}
\label{sec:iso}

We now consider the case when admissible functions $y$,
besides satisfying some boundary conditions
$y(a)=y_a$ and $y(b)=y_b$, must also satisfy
a given integral constraint. This case is known
in the literature as the isoperimetric problem
\cite{Almeida1,AlmeidaNabla}. Let
$$
\Theta(y)= \int_a^b \theta\left(x,y(x),\frac{\Box y}{\Box x}(x)\right)\, dx
= \int_a^b \theta[y](x)\, dx
$$
be a functional defined on $H^\alpha(I,\mathbb{R})$,
where $(x,y,v) \mapsto \theta(x,y,v)$ is a $C^1$ function.

\begin{Definition}
Let $y\in \partial H^\alpha(I,\mathbb{R})$. We say that $y$ is a quantum extremal for
$\Phi$ given by \eqref{funct} subject to the isoperimetric constraint
$\Theta(y) = C$, $C\in \mathbb{C}$, if for all
$n \in \mathbb{N}$ and all variations $\overline{y}=y
+\sum_{k=1}^n \varepsilon_k w_k$, where $(w_k)_{1\leq k \leq n} \in H^\beta(I,\mathbb{R})$
and $||\varepsilon|| = ||(\varepsilon_1,\ldots,\varepsilon_n)||$ is some small real number,
the condition $\Theta(\overline y)=C$ implies that
$\frac{d}{d\varepsilon_k}\left.\Phi(y + \varepsilon_k w_k)\right|_{\varepsilon = 0}=0$
for all $k\in\{1,\ldots,n\}$.
\end{Definition}

\begin{Theorem}
\label{thm:EL:ISO:P}
Let $y$ be a quantum extremal for $\Phi$ given by \eqref{funct}
subject to the isoperimetric constraint $\Theta(y) = C$, $C\in \mathbb{C}$.
If $y$ is not a quantum extremal for $\Theta$,
then there exists some $\lambda \in \mathbb{C}$ such that
$$
\frac{\partial (L-\lambda \theta)}{\partial y}[y]
-\frac{\Box}{\Box x}\left(\frac{\partial (L-\lambda \theta)}{\partial v}[y]\right)=0.
$$
\end{Theorem}

\begin{proof}
Let $\varepsilon_1,\varepsilon_2\in B_r(0,0)$, with $r$ a sufficiently small positive number,
and let $\overline y =y +\varepsilon_1 w_1 + \varepsilon_2 w_2$ be a variation of $y$
such that $w_j\in H^\beta(I,\mathbb{R})$ and $w_j(a)=w_j(b)=0$,
$j=1,2$. Define on the set $B_r(0,0)$ the functions given by
$$
\overline \Phi (\varepsilon_1,\varepsilon_2)
=\Phi(\overline y) \quad \mbox{and}\quad
\overline \Theta(\varepsilon_1,\varepsilon_2)=\Theta(\overline y)-C.
$$
Using integration by parts, one has
$$
\left.\frac{\partial \overline\Theta}{\partial \varepsilon_1}\right|_{(0,0)}
=\int_a^b \left( \frac{\partial \theta}{\partial y}[y](x)
-\frac{\Box}{\Box x} \left(\frac{\partial \theta}{\partial v}[y]\right)(x) \right) \, w_1(x) \, dx.
$$
Since $y$ is not a quantum extremal for $\Theta$, there exists a function $w_1$ for which
$(\partial \overline\Theta/\partial \varepsilon_1)|_{(0,0)}\not=0$. Thus,
we are in conditions to apply the implicit function theorem,
\textrm{i.e.}, we can write $\varepsilon_1$ as a function of $\varepsilon_2$
in a neighborhood of zero, $\varepsilon_1=\varepsilon_1(\varepsilon_2)$, such that
$$
\overline{\Theta}(\varepsilon_1(\varepsilon_2),\varepsilon_2)=0.
$$
We now adapt the Lagrange multiplier method to our complex valued functionals.
Since $\overline{\Theta}(\varepsilon_1(\varepsilon_2),\varepsilon_2)=0$
for any $\varepsilon_2$, it follows that
$$
0=\frac{d}{d\varepsilon_2}\overline{\Theta} (\varepsilon_1(\varepsilon_2),\varepsilon_2)
= \frac{d \varepsilon_1}{d \varepsilon_2} \cdot \frac{\partial \overline{\Theta}}{\partial \varepsilon_1}
+\frac{\partial\overline{\Theta}}{\partial \varepsilon_2}
$$
and so
$$
\left.\frac{d \varepsilon_1}{d \varepsilon_2}\right|_0
=- \frac{\left.\frac{\partial\overline{\Theta}}{\partial \varepsilon_2}\right|_{(0,0)}}{
\left.\frac{\partial\overline{\Theta}}{\partial \varepsilon_1}\right|_{(0,0)}}.
$$
On the other hand, since $y$ is a quantum extremal for $\Phi$ subject to the constraint $\Theta\equiv C$,
\begin{equation*}
\begin{split}
\frac{d}{d\varepsilon_2}\left.\overline{\Phi}(\varepsilon_1(\varepsilon_2),\varepsilon_2)\right|_0
&=\left.\frac{d\varepsilon_1}{d\varepsilon_2}\right|_0 \left.\frac{\partial
\overline \Phi}{\partial \varepsilon_1}\right|_{(0,0)}
+\left.\frac{\partial \overline \Phi}{\partial \varepsilon_2}\right|_{(0,0)}\\
&=\left.\frac{d\varepsilon_1}{d\varepsilon_2}\right|_0 \int_a^b\left[ \frac{\partial L}{\partial y}[y](x)
-\frac{\Box}{\Box x}\left( \frac{\partial L}{\partial v}[y]\right)(x) \right]\, w_1(x) \, dx\\
&\qquad +  \int_a^b\left[ \frac{\partial L}{\partial y}[y](x)
-\frac{\Box}{\Box x}\left( \frac{\partial L}{\partial v}[y]\right)(x)\right]\, w_2(x) \, dx \\
&=0.
\end{split}
\end{equation*}
We conclude that
$$
\left|
\begin{array}{cc}
\displaystyle\left.\frac{\partial \overline{\Phi}}{\partial \varepsilon_1}\right|_{(0,0)}
& \left.\displaystyle\frac{\partial \overline{\Theta}}{\partial \varepsilon_1}\right|_{(0,0)}\\
\displaystyle \left. \frac{\partial \overline{\Phi}}{\partial \varepsilon_2}\right|_{(0,0)}
&\displaystyle\left.\frac{\partial \overline{\Theta}}{\partial \varepsilon_2}\right|_{(0,0)}
\end{array}\right|=0.
$$
Since $(\partial \overline\Theta/\partial \varepsilon_1)|_{(0,0)}\not=0$,
there exists some $\lambda\in \mathbb{C}$ such that
$$
\left( \left.\frac{\partial \overline{\Phi}}{\partial \varepsilon_1}\right|_{(0,0)},
\left.\frac{\partial \overline{\Phi}}{\partial \varepsilon_2}\right|_{(0,0)}\right)
= \lambda \left(\left.\frac{\partial\overline{\Theta}}{\partial \varepsilon_1}\right|_{(0,0)},
\left.\frac{\partial \overline{\Theta}}{\partial\varepsilon_2}\right|_{(0,0)}\right).
$$
Consequently,
\begin{equation}
\label{equation}
\left.\frac{\partial}{\partial \varepsilon_2}(\overline \Phi
- \lambda \overline{\Theta})\right|_{(0,0)}=0.
\end{equation}
We prove the result doing the computations in equation \eqref{equation}
and using the fundamental lemma of the calculus of variations
(see, \textrm{e.g.}, \cite{Gelfand}).
\end{proof}

\begin{Remark}
In Theorem~\ref{thm:EL:ISO:P}
we are assuming that
$$
\frac{\partial \theta}{\partial v}[y](x)\cdot w(x)
$$
satisfies \eqref{nec_condition}
for all $w \in H^\beta(I, \mathbb{R})$.
\end{Remark}


\section{Dependence on a parameter}
\label{sec:parameter}

Let
\begin{equation}
\label{funct2}
\begin{array}{lcll}
\Phi: & \partial H^\alpha(I,\mathbb{R})\times \mathbb{C} & \to & \mathbb{C}\\
& (y,\xi) & \mapsto & \displaystyle \int_a^b L\left(x,y(x),\frac{\Box y}{\Box x}(x),\xi\right)\, dx.
\end{array}
\end{equation}
Introducing the operator $[\cdot]_\xi$ by
$$
[y]_\xi(x) =\left(x,y(x),\frac{\Box y}{\Box x}(x),\xi\right),
$$
we write $\Phi(y,\xi) = \int_a^b L[y]_\xi(x)\, dx$.
A pair $(y,\xi)$ is called a quantum extremal on
$H^\beta(I, \mathbb{R})\times\mathbb{C}$ for $\Phi$ as in \eqref{funct2},
if for any pair $(w,\delta)\in H^\beta(I, \mathbb{R})\times\mathbb{C}$ one has
$$
\frac{d}{d\varepsilon}\left.\Phi(y + \varepsilon w, \xi + \varepsilon \delta)\right|_{\varepsilon = 0}=0.
$$

\begin{Theorem}
\label{thm:EL:compParameter}
Let $y \in \partial H^\alpha(I,\mathbb{R})$
be such that $\Box y/\Box x \in H^\alpha(I,\mathbb{R})$.
The pair $(y,\xi)$ is a quantum extremal for $\Phi$
on $H^\beta(I, \mathbb{R})\times\mathbb{C}$ if and only if
\begin{equation}
\label{ELequation2}
\frac{\partial L}{\partial y}[y]_\xi(x)
-\frac{\Box}{\Box x}\left(\frac{\partial L}{\partial v}[y]_\xi\right)(x)=0
\end{equation}
and
\begin{equation}
\label{parameter}
\int_a^b \frac{\partial L}{\partial \xi}[y]_\xi(x)\, dx=0.
\end{equation}
\end{Theorem}

\begin{proof}
Let $(y,\xi)$ be a quantum extremal for $\Phi$.
For $(w,\delta)\in H^\beta(I, \mathbb{R})\times\mathbb{C}$ we have
\begin{equation*}
\begin{split}
0 &=  \frac{d}{d\varepsilon}\left.\Phi(y + \varepsilon w, \xi
+ \varepsilon \delta)\right|_{\varepsilon = 0}\\
&=\int_a^b\left[ \frac{\partial L}{\partial y}[y]_\xi(x)
-\frac{\Box}{\Box x}\left(\frac{\partial L}{\partial v}[y]_\xi(x)\right) \right] \cdot w \, dx
+ \int_a^b \frac{\partial L}{\partial \xi}[y]_\xi(x)\cdot \delta \, dx.
\end{split}
\end{equation*}
Taking $\delta=0$ we obtain equation \eqref{ELequation2};
taking $w \equiv 0$ we obtain \eqref{parameter}.
\end{proof}


\section{The higher-order quantum Euler--Lagrange equation}
\label{sec:higherEL}

Let $n\in\mathbb{N}$ and
$$
[y]^n(x)=\left(x,y(x),\frac{\Box y}{\Box x}(x),\ldots,\frac{\Box^n y}{\Box x^n}(x)\right).
$$
We now consider functionals of type
\begin{equation}
\label{eq:ho:functional}
\Phi(y)=\int_a^b L[y]^n(x)\, dx
\end{equation}
defined on $y \in \partial H^\alpha(I,\mathbb{R})$ with
$\Box^k y/\Box x^k$ still on $H^\alpha(I,\mathbb{R})$
for all $k \in \{1,\ldots,n \}$. We assume that the Lagrangian
$$
L=L(x,y,v_1,v_2,\ldots,v_n):[a,b]\times\mathbb{R}\times\mathbb{C}^n\to\mathbb{C}
$$
is of class $C^1$. We note that if $y\in C^n(I,\mathbb{R})$,
then \eqref{eq:ho:functional} reduces to the standard
higher-order variational functional
$$
\Phi(y)=\int_a^b L\left(x,y(x),y'(x),\ldots,y^{(n)}(x)\right)\, dx.
$$
The advantage of our quantum approach is that in order to obtain the
corresponding Euler--Lagrange equation we do not need to assume,
as it happens in the classical case \cite{Gelfand}, that
admissible functions $y$ are of class $C^{2n}$,
being necessary only continuity. As variations of $y$
we consider functions $y + \varepsilon w$ with
a small real $\varepsilon$ and
$w\in H^\beta(I,\mathbb{R})$ with $w(a)=w(b)=0$,
$\Box^k w/\Box x^k \in H^\beta(I,\mathbb{R})$,
and $\Box^k w/\Box x^k(a)=\Box^k w/\Box x^k(b)=0$
for all $k \in \{1,\ldots,n-1\}$.

\begin{Theorem}[The higher-order  quantum Euler--Lagrange equation]
\label{higherEL}
Let $\alpha, \beta \in(0,1)$ be such that $\alpha+\beta>1$ and $\beta \geq \alpha$.
A function $y \in \partial H^\alpha(I,\mathbb{R})$
with $\Box^k y/\Box x^k$ still on $H^\alpha(I,\mathbb{R})$
for all $k \in \{1,\ldots,n \}$ is a quantum extremal
of $\Phi$ on $H^\beta(I, \mathbb{R})$ if and only if
$$
\frac{\partial L}{\partial y}[y]^n(x) +\sum_{i = 1}^{n}
(-1)^i\frac{\Box^i}{\Box x^i}\left(\frac{\partial L}{\partial v_i}[y]^n\right)(x)=0
$$
for all $x\in [a,b]$.
\end{Theorem}

\begin{proof}
Using Theorem~\ref{LeibnizRule} consecutively, we deduce that
\begin{equation*}
\begin{split}
0 & = \frac{d}{d\varepsilon}\left.\Phi(y + \varepsilon w)\right|_{\varepsilon = 0}\\
&=\int_a^b \left[ \frac{\partial L}{\partial y}[y]^n(x)\cdot w(x)
+ \sum_{i = 1}^{n} \frac{\partial L}{\partial v_i}[y]^n(x)\cdot \frac{\Box^i w}{\Box x^i}(x)\right] dx\\
&= \int_a^b \left[ \frac{\partial L}{\partial y}[y]^n(x)\cdot w(x)
+ \sum_{i = 1}^{n} \frac{\partial L}{\partial v_i}[y]^n(x)\cdot
\frac{\Box}{\Box x}\left(\frac{\Box^{i-1} w}{\Box x^{i-1}}\right)(x)\right] dx\\
&= \int_a^b \left[ \frac{\partial L}{\partial y}[y]^n(x) \cdot w(x)
- \sum_{i=1}^{n} \frac{\Box}{\Box x}\left(\frac{\partial L}{\partial v_i}[y]^n\right)(x)
\cdot \frac{\Box^{i-1} w}{\Box x^{i-1}}(x)\right] dx\\
&= \int_a^b \left[ \frac{\partial L}{\partial y}[y]^n(x)
+ \sum_{i=1}^{n} (-1)^i\frac{\Box^i}{\Box x^i}\left(\frac{\partial
L}{\partial v_n}[y]^n\right)(x)\right]\cdot w(x) \, dx.
\end{split}
\end{equation*}
The result follows from the fundamental lemma of the calculus of variations.
\end{proof}

\begin{Remark}
For $n = 1$ Theorem~\ref{higherEL} reduces to Theorem~\ref{ELtheorem1}.
\end{Remark}

\begin{Remark}
In Theorem~\ref{higherEL} we are assuming that for all
$i \in \{1,\ldots,n\}$, and for all $k \in \{0,\ldots,i-1\}$,
$$
\frac{\Box^k}{\Box x^k}\left( \frac{\partial L}{\partial v_i}[y]^n\right)(x)
\cdot \frac{\Box^{i-k-1}}{\Box x^{i-k-1}} w(x)
$$
satisfies  \eqref{nec_condition} for all $w \in H^\beta(I, \mathbb{R})$.
\end{Remark}

We can easily include the case when the Lagrangian depends
on a complex parameter $\xi$.

\begin{Theorem}
Let
$$
\Phi(y,\xi)=\int_a^b L[y]_{\xi}^{n}(x)\, dx,
$$
where
$$
[y]_{\xi}^{n}(x)=\left(x,y(x),\frac{\Box y}{\Box x}(x),
\ldots,\frac{\Box^n y}{\Box x^n}(x),\xi\right),
$$
$\alpha, \beta \in(0,1)$ be such that $\alpha+\beta>1$ and $\beta \geq \alpha$,
and $y \in \partial H^\alpha(I,\mathbb{R})$
be such that $\Box^k y/\Box x^k \in H^\alpha(I,\mathbb{R})$
for $k = 1, \ldots, n$.
The pair $(y,\xi)$ is a quantum extremal for $\Phi$
on $H^\beta(I, \mathbb{R})\times\mathbb{C}$ if and only if
$$
\frac{\partial L}{\partial y}[y]_{\xi}^{n}
+\sum_{i=1}^{n} (-1)^i\frac{\Box^i}{\Box x^i}\left(\frac{\partial L}{\partial v_i}[y]_{\xi}^{n} \right)=0
$$
and
$$
\int_a^b \frac{\partial L}{\partial \xi}[y]_{\xi}^{n}(x)\, dx=0.
$$
\end{Theorem}


\section{The quantum Green theorem}
\label{sec:green}

We now prove a version of Green's theorem in the scale derivative context.
Let $\Omega\subseteq\mathbb{R}^2$ be a sufficient large set,
$R=[a,b]\times[c,d]\subseteq\Omega$ and $f\in C^0(\Omega,\mathbb{R})$.
Let $\partial R$ denote the boundary of $R$, oriented in the counterclockwise direction.

\begin{Theorem}[The quantum Green theorem]
\label{thm:HQGT}
Let $f$ and $g$ be two continuous functions
whose domains contain $\Omega$, and such that for all
$(x_1,x_2)\in[a,b]\times[c,d]$, $f(x_1,\cdot)$
and $g(\cdot,x_2)$ satisfy condition \eqref{nec_condition}.
Then,
\begin{equation}
\label{greentheorem}
\int_{\partial R}(f(x_1,x_2)dx_1+g(x_1,x_2)dx_2)
=\int\int_R\left( \frac{\Box g}{\Box x_1}(x_1,x_2)
-\frac{\Box f}{\Box x_2}(x_1,x_2) \right)dx_1dx_2.
\end{equation}
\end{Theorem}

\begin{proof}
By Theorem~\ref{Barrow}, we get
\begin{equation*}
\begin{split}
\int_{\partial R}(fdx_1+gdx_2)
&=\int_a^b f(x_1,c)-f(x_1,d)dx_1+\int_c^d g(b,x_2)-g(a,x_2)dx_2\\
&=-\int_a^b \int_c^d \frac{\Box f}{\Box x_2}(x_1,x_2)dx_2dx_1
+\int_c^d \int_a^b \frac{\Box g}{\Box x_1}(x_1,x_2)dx_1dx_2\\
&=\int\int_R\left( \frac{\Box g}{\Box x_1}(x_1,x_2)
-\frac{\Box f}{\Box x_2}(x_1,x_2) \right)dx_1dx_2.
\end{split}
\end{equation*}
The equality \eqref{greentheorem} is proved.
\end{proof}

Assume now that $f=-Fw$ and $g=Gw$,
for some continuous functions $F,G$
and $w$ such that $F,G \in H^\alpha(\Omega,\mathbb{R})$
and $w \in H^\beta(\Omega,\mathbb{R})$
with $\alpha+\beta>1$. Using Theorem~\ref{thm:HQGT}
and Theorem~\ref{LeibnizRule},
we get an integration by parts type formula:
\begin{multline}
\label{eq:intParts2d}
\int\int_R\left(G \cdot \frac{\Box w}{\Box x_1}
+F \cdot \frac{\Box w}{\Box x_2}\right)dx_1dx_2\\
=\int_{\partial R}(-Fw \, dx_1+Gw \, dx_2)
-\int\int_R\left( \frac{\Box G}{\Box x_1}
+\frac{\Box F}{\Box x_2}\right)\cdot w \, dx_1dx_2.
\end{multline}
In addition, if $w\equiv0$ on $\partial R$, then
\begin{equation}
\label{eq:intParts2d:w0}
\int\int_R\left(G \cdot \frac{\Box w}{\Box x_1}
+F \cdot \frac{\Box w}{\Box x_2}\right)dx_1dx_2
=-\int\int_R\left( \frac{\Box G}{\Box x_1}
+\frac{\Box F}{\Box x_2}\right)\cdot w \, dx_1dx_2.
\end{equation}


\section{Quantum extremals for functionals with two independent variables}
\label{sec:optimization2}

We consider functionals of type
\begin{equation}
\label{eq:2dFunct}
\begin{array}{lcll}
\Phi: & \partial H^\alpha(\Omega,\mathbb{R}) & \to & \mathbb{C}\\
& y & \mapsto & \displaystyle \int\int_R L\left(x_1,x_2,y(x_1,x_2),
\frac{\Box y}{\Box x_1}(x_1,x_2),\frac{\Box y}{\Box x_2}(x_1,x_2)\right)\, dx_1dx_2,
\end{array}
\end{equation}
where $L=L(x_1,x_2,y,v_1,v_2): R \times \mathbb{R}
\times \mathbb{C}^2 \to \mathbb{C}$ is of class $C^1$ and
$$
\partial H^\alpha(\Omega,\mathbb{R})
= \{  y \in H^\alpha(\Omega,\mathbb{R})
\, | \, y(\partial R) \mbox{ takes fixed given values}\}.
$$

\begin{Remark}
When we only consider differentiable functions $y$,
the functional \eqref{eq:2dFunct} reduces to
$$
\Phi(y)=\displaystyle \int\int_RL\left(x_1,x_2,y(x_1,x_2),
\frac{\partial y}{\partial x_1}(x_1,x_2),
\frac{\partial y}{\partial x_2}(x_1,x_2)\right)\, dx_1dx_2.
$$
\end{Remark}

As variations of $y$, we consider those functions of form $y+ \varepsilon w$
with $w \in H^\beta(\Omega, \mathbb{R})$,
$\alpha+\beta>1$ and $\beta \geq \alpha$, and $\varepsilon$ a small real number.

\begin{Definition}
We say that $y$ is a quantum extremal of $\Phi$ given by \eqref{eq:2dFunct}
on $H^\beta(\Omega, \mathbb{R})$, if for any
$w \in H^\beta(\Omega, \mathbb{R})$ one has
$$
\frac{d}{d\varepsilon}\left.\Phi(y + \varepsilon w)\right|_{\varepsilon = 0}=0.
$$
\end{Definition}

We now present the scale Euler--Lagrange equation for multiple integrals.
It has some similarities with the standard one (see, \textrm{e.g.}, \cite{Gelfand}),
replacing the partial derivatives $\partial/\partial x_j$
by the  scale derivative $\Box/\Box x_j$.

\begin{Theorem}[The  quantum Euler--Lagrange equation for double integrals]
\label{ELtheorem}
Let $\alpha,\beta\in(0,1)$ with $\alpha+\beta >1$
and $\beta \geq \alpha$, and $y \in \partial H^\alpha(\Omega,\mathbb{R})$
be such that $\Box y/\Box x_j \in H^\alpha(\Omega,\mathbb{R})$, $j \in \{1,2\}$.
The function $y$ is a quantum extremal of $\Phi$ on $H^\beta(\Omega, \mathbb{R})$
if and only if
\begin{equation}
\label{ELequation}
\frac{\partial L}{\partial y}[y](x_1,x_2)
-\frac{\Box}{\Box x_1}\left(\frac{\partial L}{\partial v_1} \circ [y]\right)(x_1,x_2)
-\frac{\Box}{\Box x_2}\left(\frac{\partial L}{\partial v_2} \circ [y]\right)(x_1,x_2) = 0
\end{equation}
for all $(x_1,x_2)\in R$, where
$$
[y](x_1,x_2)=\left(x_1,x_2,y(x_1,x_2),\frac{\Box y}{\Box x_1}(x_1,x_2),
\frac{\Box y}{\Box x_2}(x_1,x_2)\right).
$$
\end{Theorem}

\begin{proof}
Let $y$ be a quantum extremal of \eqref{eq:2dFunct}
and consider a variation $y+ \varepsilon w$ with $w\in H^\beta(\Omega, \mathbb{R})$.
Because of the given constraints on the boundary, we can only consider those variations
for which $w(\partial R)\equiv0$. Using the integration
by parts formula \eqref{eq:intParts2d:w0} we get:
\begin{equation*}
\begin{split}
0 &= \frac{d}{d\varepsilon}\left.\Phi(y + \varepsilon w)\right|_{\varepsilon = 0}\\
&= \int\int_R \left[ \frac{\partial L}{\partial y}[y](x_1,x_2)\cdot w(x_1,x_2)
+\sum_{i=1}^{2} \left(\frac{\partial L}{\partial v_i}[y](x_1,x_2)
\cdot \frac{\Box w}{\Box x_i}(x_1,x_2)\right)\right] dx_1 dx_2\\
&= \int\int_R \left[ \frac{\partial L}{\partial y}[y](x_1,x_2)
-\sum_{i=1}^{2} \frac{\Box}{\Box x_i}\left(\frac{\partial L}{\partial v_i} \circ [y]
\right)(x_1,x_2)\right] \cdot w(x_1,x_2) \, dx_1 dx_2.
\end{split}
\end{equation*}
Equation \eqref{ELequation} follows from the fundamental lemma
of the calculus of variations.
\end{proof}

We now consider the situation without the constraints on $\partial R$. Let
$$
\begin{array}{lcll}
\Psi: & H^\alpha(\Omega,\mathbb{R}) & \to & \mathbb{C}\\
    & y & \mapsto & \Phi(y).\\
\end{array}
$$

\begin{Theorem}[The quantum Euler--Lagrange equation for double integrals
and natural boundary conditions]
\label{ELtheorem:di}
Let $\alpha,\beta\in(0,1)$ with $\alpha+\beta >1$ and $\beta \geq \alpha$,
and $y \in H^\alpha(\Omega,\mathbb{R})$ be such that
$\Box y/\Box x_j \in H^\alpha(\Omega,\mathbb{R})$, $j \in \{1,2\}$.
Then, $y$ is a quantum extremal of $\Psi$ on $H^\beta(\Omega,\mathbb{R})$
if and only if the following conditions hold:
\begin{enumerate}
\item $\displaystyle\frac{\partial L}{\partial y}[y](x_1,x_2)
-\frac{\Box}{\Box x_1}\left(\frac{\partial L}{\partial v_1} \circ [y]\right)(x_1,x_2)
-\frac{\Box}{\Box x_2}\left(\frac{\partial L}{\partial v_2} \circ [y]\right)(x_1,x_2)
=0$  $\forall\, (x_1,x_2) \in R$;
\item $\displaystyle\frac{\partial L}{\partial v_1}\left( a,x_2,y(a,x_2),
\frac{\Box y}{\Box x_1}(a,x_2),\frac{\Box y}{\Box x_2}(a,x_2)\right)
=0$ $\forall\, x_2\in[c,d]$;
\item $\displaystyle\frac{\partial L}{\partial v_1}\left( b,x_2,y(b,x_2),
\frac{\Box y}{\Box x_1}(b,x_2),\frac{\Box y}{\Box x_2}(b,x_2) \right)
=0$ $\forall\, x_2\in[c,d]$;
\item $\displaystyle\frac{\partial L}{\partial v_2}\left( x_1,c,y(x_1,c),
\frac{\Box y}{\Box x_1}(x_1,c),\frac{\Box y}{\Box x_2}(x_1,c) \right)
=0$ $\forall\, x_1\in[a,b]$;
\item $\displaystyle\frac{\partial L}{\partial v_2}\left( x_1,d,y(x_1,d),
\frac{\Box y}{\Box x_1}(x_1,d),\frac{\Box y}{\Box x_2}(x_1,d) \right)
=0$ $\forall\, x_1\in[a,b]$.
\end{enumerate}
\end{Theorem}

\begin{proof}
Similarly to what was done in the proof of Theorem~\ref{ELtheorem},
but using \eqref{eq:intParts2d} instead of \eqref{eq:intParts2d:w0},
we deduce that
\begin{multline}
\label{proof}
0=\int\int_R \left[ \frac{\partial L}{\partial y}
-\frac{\Box}{\Box x_1}\left(\frac{\partial L}{\partial v_1} \right)
-\frac{\Box}{\Box x_2}\left(\frac{\partial L}{\partial v_2}\right) \right]
\cdot w\, dx_1 dx_2\\
+ \int_{\partial R} \left[ -\frac{\partial L}{\partial v_2}
\cdot w \, dx_1 + \frac{\partial L}{\partial v_1}\cdot w \, dx_2 \right].
\end{multline}
In particular, we may consider those $w$ that are zero on $\partial R$,
and we obtain equality \eqref{ELequation} of item \textit{1}.
Replacing \eqref{ELequation} into equation \eqref{proof}, it follows that
$$
\begin{array}{ll}
0 & =\displaystyle \int_{\partial R} \left[
-\frac{\partial L}{\partial v_2}\cdot w \, dx_1
+ \frac{\partial L}{\partial v_1}\cdot w \, dx_2\right] \\
& = -\displaystyle \int_a^b \frac{\partial L}{\partial v_2}\left( x_1,c,y(x_1,c),
\frac{\Box y}{\Box x_1}(x_1,c),\frac{\Box y}{\Box x_2}(x_1,c) \right)\cdot w(x_1,c) \, dx_1\\
& \quad +\displaystyle \int_a^b \frac{\partial L}{\partial v_2}\left( x_1,d,y(x_1,d),
\frac{\Box y}{\Box x_1}(x_1,d),\frac{\Box y}{\Box x_2}(x_1,d) \right)\cdot w(x_1,d) \, dx_1\\
& \quad+\displaystyle \int_c^d  \frac{\partial L}{\partial v_1}\left( b,x_2,y(b,x_2),
\frac{\Box y}{\Box x_1}(b,x_2),\frac{\Box y}{\Box x_2}(b,x_2) \right)  \cdot w(b,x_2) \, dx_2\\
& \quad-\displaystyle \int_c^d   \frac{\partial L}{\partial v_1}\left( a,x_2,y(a,x_2),
\frac{\Box y}{\Box x_1}(a,x_2),\frac{\Box y}{\Box x_2}(a,x_2) \right) \cdot w(a,x_2) \, dx_2.
\end{array}
$$
The result is proved with appropriate choices of $w$.
\end{proof}


\section{An application}
\label{sec:example}

A membrane is a portion of surface, plane at rest, with potential energy proportional
to the change in area. Suppose that the membrane at rest covers an area $R\subseteq \mathbb{R}^2$,
and in the boundary the membrane is fixed. If $u(\cdot, \cdot)$ denotes the deformation normal
to the equilibrium, then the required potential energy is given by the expression
\begin{equation}
\label{eq:funct:appl}
\frac12 \int \int _R \left(\frac{\partial u}{\partial x_1}\right)^2
+\left(\frac{\partial u}{\partial x_2}\right)^2dx_1dx_2.
\end{equation}
The equilibrium is attained at the Euler--Lagrange extremal of this double integral.
For more on the subject, we refer the reader to \cite{Courant}.

The functional \eqref{eq:funct:appl} is considered for
continuous functions $u$ that possess continuous first and second order derivatives.
These can be very restrictive assumptions.
Indeed, it is well known that often real phenomenons are nondifferentiable,
and so not suitable for treatment under the standard variational approach.
Assume that we wish to include in the above problem nondifferentiable functions,
where the velocity is described by a  quantum derivative. In the  scale calculus context,
we consider the functional
$$
\Phi(u)=\frac12\int\int_R \left(\frac{\Box u}{\Box x_1}\right)^2
+\left(\frac{\Box u}{\Box x_2}\right)^2 dx_1dx_2
$$
defined on $\partial H^\alpha(\Omega,\mathbb{R})$.
Then, by Theorem~\ref{ELtheorem},
any quantum extremal $u$ must fulfill the equation
$$
\frac{\Box}{\Box x_1}\left(\frac{\Box u}{\Box x_1} \right)
+\frac{\Box}{\Box x_2}\left(\frac{\Box u}{\Box x_2} \right)=0.
$$


\section*{Acknowledgments}

Work supported by the \emph{Portuguese Foundation
for Science and Technology} (FCT) through the R\&D Unit
\emph{Center for Research and Development
in Mathematics and Applications} (CIDMA).



\end{document}